\documentclass[twocolumn]{autart}
\usepackage{epsfig}
\usepackage{amssymb,epsfig,color,cite,amsmath,amsfonts,mathrsfs,algorithm,algorithmic}

\newtheorem{theorem}{Theorem}
\newtheorem{lemma}{Lemma}
\newtheorem{remark}{Remark}

\begin{document}

\begin{frontmatter}

\title{Distributed Observers Design for Leader-Following Control of
Multi-Agent Networks {(Extended Version)}}
\date{}
\footnotetext{This work was supported in part by the NNSF of China
under Grants 60425307, 60628302, 10472129 and 60221301, and in
part by the US NSF Grant No. ECS-0322618 and Grant No.
ECS-0621605}

\author[beijing]{Yiguang Hong}\ead{yghong@iss.ac.cn},    
\author[hongkong]{Guanrong Chen}\ead{eegchen@cityu.edu.hk},               
\author[usa]{Linda Bushnell}\ead{lb2@u.washington.edu}  

\address[beijing]{Institute of Systems
Science, Chinese Academy of Sciences, Beijing 100080, China }
\address[hongkong]{Department
of Electronic Engineering, City University of Hong Kong, China}
\address[usa] {Department of Electrical
Engineering University of Washington, Seattle, WA 98195 USA}


\begin{abstract}
This paper is concerned with a leader-follower problem for a
multi-agent system with a switching interconnection topology.
Distributed observers are designed for the second-order
follower-agents, under the common assumption that the velocity of
the active leader cannot be measured in real time. Some dynamic
neighbor-based rules, consisting of distributed controllers and
observers for the autonomous agents, are developed to keep
updating the information of the leader. With the help of an
explicitly constructed common Lyapunov function (CLF), it is
proved that each agent can follow the active leader. Moreover, the
tracking error is estimated even in a noisy environment.  Finally,
a numerical example is given for illustration.
\end{abstract}

\end{frontmatter}

\section{Introduction}

Collective behaviors of large numbers of autonomous individuals
have been extensively studied from different points of view.  A
multi-agent network provides an excellent model for describing and
analyzing complex interconnecting behaviors, with applications in
many disciplines of physics, biology, and engineering (Okubo,
1986; Kang, Xi, \& Sparks, 2000; Lin, Broucke, \& Francis, 2004;
Ren \& Beard, 2005). Many interesting agent-related problems are
under investigation and leader-following is one of the main
research topics (Olfati-Saber, 2006; Shi, Wang, \& Chu, 2006;
Hong, {\it et al}, 2007). Neighbor-based rules are widely applied
in multi-agent coordination, inspired originally by the
aggregations of groups of individual agents in nature. In
practice, multi-agent systems typically need distributed sensing
and control due to the constraints on, or the confluence of
actuation, communication and measurement.

Distributed estimation via observers design for multi-agent
coordination is an important topic in the study of multi-agent
networks, with wide applications especially in sensor networks and
robot networks, among many others. Yet, very few theoretic results
have been obtained to date on distributed observers design and
measurement-based dynamic neighbor-based control design.
Nevertheless, one may find in the literature that Fax and Murray
(2004) reported some results concerning with distributed dynamic
feedback of special multi-agent networks, and Hong {\em et al.}
(2006) proposed an algorithm for distributed estimation of the
active leader's unmeasurable state variables, to name just a
couple.

The motivation of this work is to expand the conventional
observers design to the distributed observers design for a
multi-agent system where an active leader to be followed moves in
an unknown velocity. The continuous-time agent models considered
here are second-order, different from those first-order ones
discussed in (Hong, Hu, \& Gao, 2006). Also, switching inter-agent
topologies are taken into account here, for which a common
Lyapunov function (CLF) will be constructed. As commonly known, it
is not an easy task to construct a CLF for a switching system,
especially when the dimension of the system is high. The approach
adopted here is to reduce the order of the distributed observers
so as to reduce the dimension of the whole system, which can
significantly simplify the construction of the needed CLF.

\section{Preliminaries}

Consider a system consisting of one leader and $n$
agent-followers. A simple and undirected graph $\mathcal{G}$
describes the network of these $n$ agents, and $\bar{\mathcal{G}}$
denotes the graph that consists of $\mathcal{ G}$ and also the
leader, where some agents in $\mathcal{ G}$ are connected to the
leader via directed edges.  The graph $\mathcal{ G}$ is allowed to
have several components, within every such component all the
agents are connected via undirected edges in some topologies. The
graph $\bar{\mathcal{G} }$ of this multi-agent system is said to
be connected if at least one agent in each component of $\mathcal{
G}$ is connected to the leader by a directed edge.

In practice the relationships among neighboring agents may vary
over time, and their interconnection topology may also be
dynamically changing. Suppose that there is an infinite sequence
of bounded, non-overlapping, continuous time-intervals
$[t_{j},t_{j+1})$, $j=0,1,\cdots$, say starting at $t_{0}=0$, over
which $\sigma_N: [0,\infty)\rightarrow \mathcal{P}=
\{1,2,\cdots,N\}$ is a piecewise constant switching signal for
each $N$, defined at successive switching times. To avoid infinite
switching during a finite time interval, assume as usual that
there is a constant $\tau$ with $t_{j+1}-t_j\geq \tau$ for all
$j\geq 0$.

Let $\mathcal{N}_{i}(t)$ be the set of labels of those agents that
are neighbors of agent $i$ at time $t$. Moreover, $a_{ij}=a_{ji}$
($i=1,...,n;\,j=1,...,n$) with $a_{ii}=0$ denote the nonzero
interconnection weights between agent $i$ and agent $j$. Then, the
Laplacian of the weighted graph $\mathcal{ G}$ is denoted by $L$
(see Godsil \& Royle (2001) for the details). Moreover let
$\mathcal{N}_{0}(t)$ be the set of labels of those agents that are
neighbors of the leader at time $t$, and the nonzero connection
weight between agent $i$ and the leader (simply labelled 0),
denoted by $b_i$ for $i=1,...,n$.

Assume that the leader is active, in the sense that its state
keeps changing throughout the entire process, with dynamics
described as follows:
\begin{equation}
\label{leader1}
\begin{cases} \dot{x}_0=v_0,\quad x_0\in R^m \\ \dot
v_0=u_0,\quad v_0\in R^m\\
y=x_0\end{cases}
\end{equation}
where $x_0$ is the position, $v_0$ is the velocity, and $y$ is the
only measurable variable. This work is to expand the conventional
observer design (where the input is somehow known) to a
neighbor-based observer design. In some practical cases, the
velocity $v_0$ is hard to measure in real time, but the input
$u_0(t)$ may be regarded as some given policy known to all the
agents.

The dynamics of follower-agent $i$ is described by
\begin{equation}
\label{agent1}\begin{cases} \dot{x}_i=v_i+\delta_i^1,\quad x_i\in R^m \\
\dot v_i=u_i+\delta_i^2,\quad v_i\in R^m,\quad
i=1,...,n\end{cases}
\end{equation}
where $\delta_i^j(t)$ $(j=1,2)$ disturbances, and $u_i$
$(i=1,...,n)$, the interaction inputs. As usual, we assume that
$|\delta_i^j| \leq\Delta<\infty$ for all $j=1,2;\,i=1,...,n$. The
problem is to let all the follower-agents keep the same pace of
the leader. Without loss of generality in the following analysis,
let $m=1$ just for notational simplicity.

The following lemma (Horn \& Johnson, 1985) will be useful later.

\begin{lemma}
\label{lem1} Consider a symmetric matrix
$$
D=\begin{pmatrix}A&E \\ E^{T}&C \end{pmatrix}
$$
where $A$ and $C$ are square. Then $D$ is positive definite if and
only if both $A$ and $C- E^{T}A^{-1}E$ are positive definite.
\end{lemma}

\section{Main Results}

Since all agents cannot obtain the value of $v_0$ of the leader in
real time, they have to estimate it throughout the process. To be
more specific, denote by $\hat v_i$ an estimate of $v_0$ by agent
$i$ ($i=1,...,n$). Then, for agent $i$ to track the active leader,
the following neighbor-based rule is proposed:
$$
u_i=u_0-k[v_i-\hat v_i]
$$
\begin{equation}
\label{contr}-l\left[\sum_{j\in \mathcal{ N}_i(t)}a_{ij}(x_i-
x_j)+\sum_{i\in \mathcal{ N}_0(t)}b_i(x_i-x_0)\right],
\end{equation}
for $i=1,\cdots,n$ and constants $k,l>0$ to be determined, along
with the following distributed ``observer''
\begin{equation}
\label{update1} \dot {\hat v}_i=u_0-\frac{l}{k^2}\left[\sum_{j\in
\mathcal{ N}_i}a_{ij}(x_i- x_j)+\sum_{i\in \mathcal{ N}_0}
b_i(x_i-x_0)\right]
\end{equation}
for $i=1,\cdots,n$. Clearly, (\ref{contr}) and (\ref{update1})
contain only local information (from the agent itself and its
neighbors).

\begin{remark}
Note that there was a typo in Hong et al. (2008), where (\ref{update1}) was written as
$$
\dot {\hat v}_i=u_0-\frac{l}{k}\left[\sum_{j\in
\mathcal{ N}_i}a_{ij}(x_i- x_j)+\sum_{i\in \mathcal{ N}_0}
b_i(x_i-x_0)\right]
$$
and then $\zeta=(\hat v_1\, \cdots \, \hat
v_n)^T-v_0\mathbf{1}$.  The typo has been corrected by some papers such as Wang \& Hong (2009). In fact, the typo can be easily found and fixed when we simply go through a reverse procedure from the key equation (\ref{model2}).
\end{remark}

\begin{remark}
Hong et al. (2008) also tended to provide a general structure for the distributed observer design.  In fact, the design (\ref{update1}) can be generalized as
$$
\dot {\hat v}_i=u_0-k_0\left[\sum_{j\in
\mathcal{ N}_i}a_{ij}(x_i- x_j)+\sum_{i\in \mathcal{ N}_0}
b_i(x_i-x_0)\right]
$$
with an adjustable parameter $k_0$ to replace $\frac{l}{k^2}$.   Then we have one more parameter for more freedom to design the distributed observer for (\ref{agent1}).
\end{remark}

In Hong et al. (2006), the ``observer" has the same
dimension as the agents in a single-integrator form. Here, both
the leader and the follower-agents are described by a double
integrator, but the ``observer" is of the first-order. In fact, it
is preferred to have a one-dimensional reduced-order ``observer"
(\ref{update1}) instead of second-order ``observers"
(corresponding to the second-order agents), regarding possible
technical difficulty in constructing a CLF for the higher-order
system later on.

At first, consider the system in a noise-free environment; that
is, $\Delta=0$ ({\it i.e.}, $\delta_i^j=0$ for all $j=1,2; \,
i=1,...,n$).

\begin{theorem}
\label{thm1} Consider the leader (\ref{leader1}) and $n$ agents
(\ref{agent1}) with $\Delta=0$. In each time interval
$[t_{i},t_{i+1})$, if the entire graph is connected, then there
are constants $k$ and $l$ such that controller (\ref{contr}) with
``observer" (\ref{update1}) together yields
\begin{equation}
\label{result11} \lim_{t\rightarrow \infty}
|x_i(t)-x_0(t)|=0,\quad \lim_{t\rightarrow \infty}
|v_i(t)-v_0(t)|=0;
\end{equation}
namely, the agents can follow the leader (in the sense of both
position and speed).
\end{theorem}

\noindent Proof:\ For simplicity, set $ \xi= (x_1\, \cdots\,
x_n)^T-x_0\mathbf{1}$, $\eta=( v_1\, \cdots \,
v_n)^T-v_0\mathbf{1}$, and $\zeta=k(\hat v_1\, \cdots \, \hat
v_n)^T-kv_0\mathbf{1}$\footnote{ Hong, et al (2008) carelessly wrote $\zeta=(\hat v_1\, \cdots \, \hat
v_n)^T-v_0\mathbf{1}$, which could not yield (\ref{model2}), so we could find the typo easily from (\ref{model2}). }, where $\mathbf{1}=(1\,\cdots\, 1)^T\in
R^n$. Then, in the case of $\Delta=0$, the closed-loop system with
(\ref{contr}) and (\ref{update1}) can be written as
$$
\begin{cases}
\dot\xi=\eta \\
\dot \eta=-l(L_{\sigma}+B_{\sigma})\xi-k\eta +\zeta, \\
\dot \zeta =- \frac{l}{k}(L_{\sigma}+B_{\sigma})\xi
\end{cases}\; z=\begin{pmatrix} \xi \\
\eta \\ \zeta\end{pmatrix}\in R^{3n}
$$
or, in a compact form,
\begin{equation}
\label{model2} \dot z=F_{\sigma} z,\quad
F_{\sigma}=\begin{pmatrix} 0&I&0\\ -lH_{\sigma} & -kI &I\\
-\frac{l}{k}H_{\sigma} &0&0\end{pmatrix},\;
H_\sigma=L_\sigma+B_\sigma
\end{equation}
where the switching signal $\sigma: [0,\infty)\rightarrow
\mathcal{P}$ is piecewise constant, $B_{\sigma}$ is an $n\times n$
diagonal matrix whose $i$th diagonal element is either $b_i$ (if
agent $i$ is connected to the leader) or 0 (if it is not
connected), and $L_{\sigma}$ is the Laplacian of the $n$ agents.
In each time interval, $L_p$ and $B_p$ are time-invariant for some
$p\in \mathcal{ P}$.

By Lemma 3 of Hong et al. (2006), $H_{p}=L_p+B_p$ is
positive definite since the switching graph $\bar {\mathcal{ G}}$
remains being connected. Moreover, once $n$ is given,
$\bar\lambda$ and $\underline\lambda$, denoting the maximum and
minimum positive eigenvalues of all the positive definite matrices
$H_p$, $p\in \mathcal{ P}$, are fixed and depend directly on the
given constants $a_{ij}$ and $b_i$, $i=1,...,n;\,j=1,...,n$.
Select
\begin{equation}
\label{l} l \geq {2}/{\underline{\lambda}},\quad k\geq 4+\bar
\lambda\, l.
\end{equation}

Here, for system (\ref{model2}), a CLF is constructed as
$V(z)=z^{T}(t)Pz(t)$, with
\begin{equation}
\label{matrixp} P=\begin{pmatrix}kI &I &- \frac{k}{2} I\\
I& I& -\frac{1}{2}I\\ - \frac{k}{2} I&-\frac{1}{2}I&
\frac{k}{2}I\end{pmatrix}
\end{equation}
which is positive definite due to (\ref{l}).

Take an interval $[t_i, t_{i+1})$ into consideration. According to
the assumed conditions, the graph associated with $H_p$ for some
fixed $p\in \mathcal{P}$ is connected and time-invariant. The
derivative of $V(z)$ is given by
\begin{equation}
\label{dotv}
\dot{V}(z)\,\Big|_{(\ref{model2})}=z^T(F_{p}^{T}P+PF_{p})z
-z^{T}Q_{p}z
\end{equation}
where
$$
Q_{p}=\begin{pmatrix}
2lH_p-l H_p & l H_{p}-\frac{l}{2k}H_p & -I\\
l H_{p}-\frac{l}{2k}H_p&2(k-1)I& -I\\ -I& -I & I
\end{pmatrix}
$$
Set
$$
K_p=\begin{pmatrix}
2(k-1)I& -I\\
-I& I
\end{pmatrix} - \begin{pmatrix}
\frac{(2-1/k)^2}{4} l H_p & \frac{2-1/k}{2}I\\
\frac{2-1/k}{2}I & \frac{1}{l}H_p^{-1}
\end{pmatrix}
$$
$$
=\begin{pmatrix}
2(k-1)I-\frac{(2-1/k)^2}{4} l H_p& -\frac{4-1/k}{2}I\\
-\frac{4-1/k}{2}I& I- \frac{1}{l}H_p^{-1}
\end{pmatrix}
$$
From (\ref{l}), one has $I- \frac{1}{l}H_p^{-1}\geq I/2$, and
again by (\ref{l}), one can see that $K_p$ is a positive definite
matrix according to Lemma \ref{lem1}. Thus, by recalling Lemma
\ref{lem1} again, one can verify the positive definiteness of
$Q_p$.

It follows that there is a constant $\beta$, independent of the
selection of the time intervals, such that $\dot V(z)\leq -2\beta
V(z)$, {\it i.e.}, $V(z(t))\leq V(z(t_i))e^{-2\beta (t-t_i)},\;
\forall\,t\in [t_i, t_{i+1})$. Consequently,
\begin{equation}
\label{expf} V(z(t))\leq V(z(0))\,e^{-2\beta t}, \quad t_0=0
\end{equation}
which implies (\ref{result11}). \hfill\rule{4pt}{8pt}

Next, return to system (\ref{agent1}) with $\Delta\neq 0$. Let $T$
be a positive constant and take a sequence of intervals
$[T_j,T_{j+1})$ with $T_0=0$ and $T_{j+1}=T_j+T$. Then we have:

\begin{theorem}
\label{thm2} In each time interval $[T_j,T_{j+1})$, if the total
period over which the entire graph is connected is sufficient
large, then there is a constant $c_\delta>0$ with $\lim_{\Delta\to
0} c_{\delta}=0$, such that
\begin{equation}
\label{result12} \lim_{t\rightarrow \infty} |x_i(t)-x_0(t)|\leq
c_\delta, \quad \lim_{t\rightarrow \infty} |v_i(t)-v_0(t)|\leq
c_\delta
\end{equation}
for the multi-agent system (\ref{leader1})-(\ref{agent1}) with
local feedback (\ref{contr}) and (\ref{update1}).  In other words,
the tracking error of each agent is bounded.
\end{theorem}

\noindent Proof:\ Following the proof of Theorem \ref{thm1}, one
can obtain
\begin{equation}
\label{unknown1}
\begin{cases}
\dot\xi=\eta +\delta_1 \\
\dot \eta=-l(L_{\sigma}+B_{\sigma})\xi-k\eta +\zeta +\delta_2 \\
\dot \zeta =- \frac{l}{k}(L_{\sigma}+B_{\sigma})\xi
\end{cases}
\end{equation}
where $\delta_j=(\delta_1^j\; \cdots\; \delta_{n-1}^j\;
\delta_n^j)^T\in R^n,\; j=1,2;$, or, in a compact form, $\dot
z=F_{\sigma}z +\delta$, where $F_{\sigma}$ was defined in
(\ref{model2}) and  $\delta= (\delta_1^T\; \delta_2^T\; 0\;
\cdots\; 0)^T\in R^{3n}$.

Still take $V(z)=z^{T}Pz$ with $P$ given in (\ref{matrixp}). Each
interval $[T_j,T_{j+1})$ may consist of a number of subintervals
(still denoted by $[t_{i},t_{i+1})$ for some $i$), during which
the graph associated with $H_p$ for some $p\in \mathcal{P}$ is
connected and unchanged. Hence, we still have positive definite
matrix $Q_p$ and a constant $\beta>0$ given in the proof of
Theorem 1. Consequently, one has
\begin{equation}
\label{vv1} V(z(t_{i+1}))\leq e^{-\beta
(t_{i+1}-t_i)}\,V(z(t_{i}))+ \beta_0\Delta^{2}
\end{equation}

On the other hand, during period $[t_{\iota},t_{\iota+1})$ for
some $\iota$, the graph associated with $H_{p'}$ for some $p'\in
\mathcal{P}$ is unconnected. So, there is a constant $\alpha>0$
such that $-z^{T} Q_{p'} z \leq \frac{\alpha}{2} V(z)$.
Consequently, there is a constant $\alpha_0>0$ such that $ \dot
V(z)\leq \alpha V(z)+ \alpha_0\Delta^{2} $. Denote by $t_{j}^{d}$
the total length of all the intervals $[t_{\iota},t_{\iota+1})$ in
$[T_j,T_{j+1})$ during which the graph is unconnected, and let
$t_d=\max_j \,\{t_{j}^{d}\}$. Then,
$$
V(z(t_{\iota+1}))\leq e^{\alpha (t_{\iota+1}-t_\iota)}\,
V(z(t_{\iota}))+ \frac{\alpha_0}{\alpha}(e^{\alpha
t_d}-1)\Delta^{2}
$$
It follows that, during the time interval $[T_j,T_{j+1})$,
\begin{equation}
V(z(T_{j+1}))\leq e^{-\beta (T-t_{j}^{d})+\alpha
t_{j}^{d}}\,V(z(T_{j}))+\bar{g}\Delta^2
\end{equation}
where $\bar{g}=\frac{e^{(m+1)t_d}-1}{e^{t_d}-1}\max\left\{\beta_0,
\frac{\alpha_0}{\alpha}(e^{\alpha t_d}-1)\right\}$. If the total
period over which the graph is connected (that is, $T-t^{d}$) is
sufficiently large, then $\varepsilon=e^{-\beta (T-t^{d})+\alpha
t^{d}}< 1$. Consequently,
$$
V(z(T_{j+1})) \leq
\varepsilon^{i+1}V(z(T_{0}))+\frac{1-\varepsilon^{i+1}}
{1-\varepsilon}\,\bar{g}\Delta^2.
$$
As $j \to\infty\; (i.e.,\; t \to \infty)$, $ V(z) \leq
\frac{1}{1-\varepsilon}\,\bar{g}\Delta^2$, which implies the
conclusion. \hfill\rule{4pt}{8pt}

Here a simulation result is presented for illustration.
Consider a multi-agent system with one leader and four followers.
The interconnection topology is time-varying of switching period
0.2 between two graphs $\bar{\mathcal{G}}_i\; (i=1,2)$ described
as follows. The Laplacians for the two subgraphs
${\mathcal{G}}_i\; (i=1,2)$ of the four followers are
$$
L_1=\begin{pmatrix} 2& -1& -1& 0\\ -1& 2& 0& -1\\-1& 0& 1& 0\\ 0&
-1& 0& 1\end{pmatrix}\quad L_2=\begin{pmatrix} 1 & -1& 0& 0\\-1&
1& 0& 0\\ 0& 0& 1& -1\\ 0& 0& -1& 1\end{pmatrix}
$$
and the diagonal matrices for the interconnection relationship
between the leader and the followers are
$$
B_1=diag\{ 1\, 0 \, 0 \,0\}\in R^{4\times 4}\quad  B_2=diag\{ 1\,
0\, 1\, 0\}\in R^{4\times 4}
$$
The numerical results are obtained with $k=200$, $l=40$, and
$u_0=\cos (t)$. Fig. 1 shows that the follower-agents can track
the leader in the noise-free case and there are some bounded
errors in the case with disturbances. This further verifies the
above analysis.

\begin{figure}[ht]
\centerline{\epsfig{figure=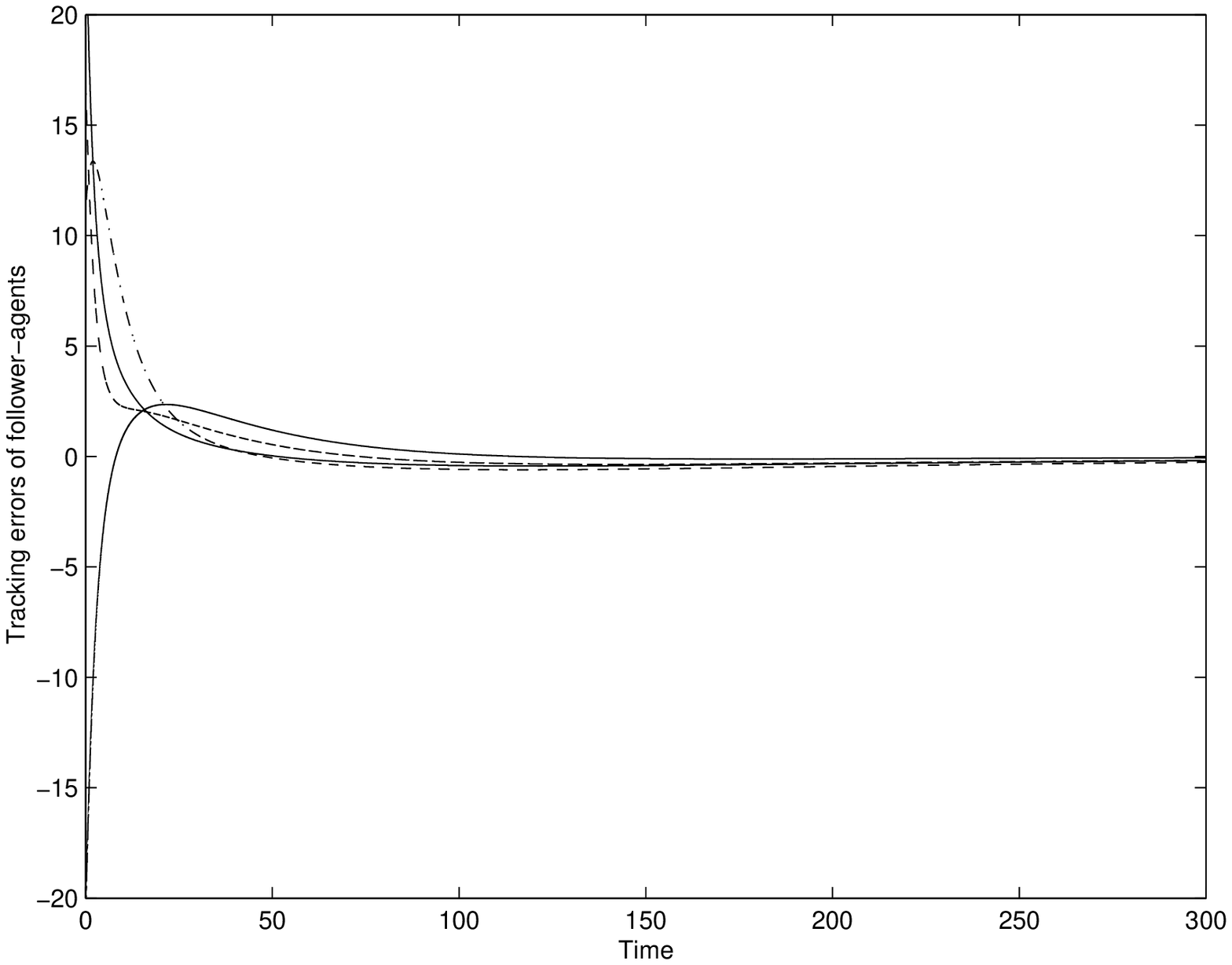,
width=0.5\linewidth=0.3}\epsfig{figure=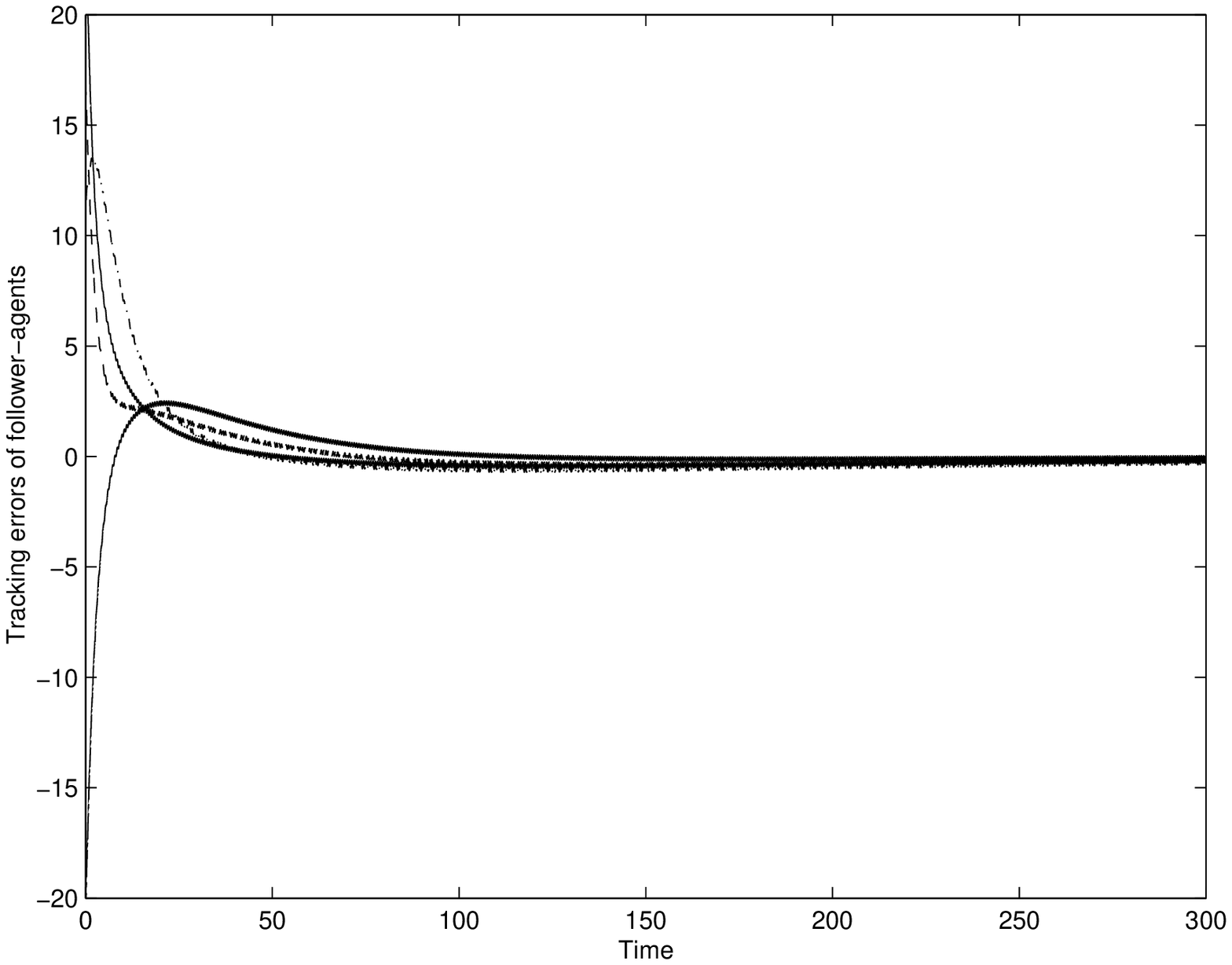,
width=0.5\linewidth=0.3}} \caption{Position tracking errors of
four followers: noise-free (on the left) and with disturbance
$\delta_i^j=\sin 50t, \, i=1,...,4, j=1,2$ (on the right)}
\end{figure}

\section{Conclusions}
This paper discussed a group of mobile agents with an active
leader moving with an unknown velocity. A neighbor-based observer
design approach was proposed, along with a dynamic coordination
rule developed for each autonomous agent. It was proved that this
distributed control guarantees the leader-following in a switching
network topology. Moreover, the tracking error has been evaluated,
even in a noisy environment.

\begin{ack}
This work was supported in part by the NNSF of China under Grants
60425307, 10472129, 50595411, and 60221301, and in part by the US
NSF Grant No. ECS-0322618 and Grant No. ECS-0621605.
\end{ack}

\section*{References}
\begin{description}
\item Fax, A., \& Murray,  R. M. (2004). Information flow and
cooperative control of vehicle formations, {\em IEEE Trans. on
Automatic Control}, 49(9), 1465-1476, .

\item Godsil C. \& Royle, G. (2001). {\it Algebraic Graph Theory,}
New York: Springer-Verlag.

\item Hong, Y., Gao, L., Cheng, D., \& Hu, J. (2007).
Lyapunov-based approach to multi-agent systems with switching
jointly-connected interconnection, {\it IEEE Trans. Automatic
Control}, 52(5), 943-948.

\item Hong, Y., Hu, J., \& Gao, L. (2006). Tracking control for
multi-agent consensus with an active leader and variable topology,
{\it Automatica}, 42(7), 1177-1182.

\item Horn, R., \& Johnson, C. (1985). {\it Matrix Analysis,} New
York: Cambbridge Univ. Press.

\item Kang, W., Xi, N., \& Sparks, A. (2000). Formation control of
autonomous agents in 3D workspace, {\em Proc. of IEEE Int. Conf.
on Robotics and Automation}, 1755-1760, San Francisco, CA.

\item Lin, Z., Broucke, M., \& Francis, B. (2004). Local control
strategies for groups of mobile autonomous agents, {\it IEEE
Trans. Automatic Control}, 49(4), 622-629.

\item Okubo, A. (1986). Dynamical aspects of animal grouping:
swarms, schools, flocks and herds, {\it Advances in Biophysics},
22, 1-94.

\item Olfati-Saber, R. (2006). Flocking for multi-agent dynamic
systems: algorithms and theory, {\em IEEE Trans. on Automatic
Control}, 51(3): 410-420.

\item Ren, W., \& Beard, R. (2005). Consensus seeking in
multi-agent systems using dynamically changing interaction
topologies, {\em IEEE Trans. Automatic Control}, 50(4), 665-671.

\item Shi, H., Wang, L., \& Chu, T. (2006). Virtual leader
approach to coordinated control of multiple mobile agents with
asymmetric interactions, {\it Physica D}, 213, 51-65.

\item Hong, Y., Chen, G., \& Bushnell, L. (2008). Distributed observers design for leader-following control of
multi-agent networks, {\it Automatica}, 44 (3), 846-850.

\item Wang, X.,\& Hong, Y. (2009). Distributed observers for tracking a moving target by cooperative multiple agents with time delays, {\it Proc. of ICCAS-SICE 2009}, 982-987, Japan.
\end{description}
\end{document}